\documentclass[a4paper,12pt]{article}
\usepackage{amsmath}
\usepackage{amsthm}
\usepackage{amssymb}
\usepackage{amsfonts}
\usepackage{amscd}
\theoremstyle{plain}
	\newtheorem{thm}{Theorem}[section]
	\newtheorem*{thm*}{Theorem}
	\newtheorem*{cor}{Corollary}
	\newtheorem{lem}[thm]{Lemma}
	\newtheorem{prop}[thm]{Proposition}
\theoremstyle{definition}
	\newtheorem{defn}{Definition}[section]
	\newtheorem{conj}{Conjecture}[section]

\theoremstyle{remark}
	\newtheorem*{rem}{Remark}

	\newtheorem*{pf}{Proof}
\def\C{{\mathbb C}}
\def\Z{{\mathbb Z}}
\def\P{{\mathbb P}}
\def\R{{\mathbb R}}
\def\Q{{\mathbb Q}}
\def\N{{\mathbb N}}
\begin{document}
\title{Primitive Forms, \\
	Topological LG Models Coupled to Gravity\\
	and \\
	Mirror Symmetry}
\author{Atsushi Takahashi
	\thanks{{\it E-mail address}~: atsushi@kurims.kyoto-u.ac.jp} \\
	\\
	Research Institute for Mathematical Sciences,\\
	Kyoto University,\\
	Kyoto 606-01, Japan}
\date{}
\maketitle
\begin{abstract}
In this paper, we will describe the mathematical foundation of
topological Landau--Ginzburg (LG) models coupled to gravity at genus $0$
in terms of primitive forms.
In particular, we can give a natural explanation for the large phase space.
We also discuss the mirror symmetry for Calabi--Yau manifolds and $\C\P^1$ in
our context.
We see that the mirror partner of $\C\P^1$ is the theory of primitive form
associated to $f=z+qz^{-1}$.
\end{abstract}
%
\section*{Introduction}
~~~
Mirror symmetry (for detail see \cite{gy:1}\cite{ya:1}) is the one of the
most important phenomena in algebraic geometry and string theory.
A pair of manifolds $(X,Y)$ is said to be mirror pair if
A-twisted sigma model coupled to gravity on $X$ is isomorphic to B-twisted
sigma model coupled to gravity on $Y$ \cite{wi:3}.
The correlation functions of A-twisted sigma model coupled to
gravity have been mathematically established as the theory of Gromov--Witten
invariants \cite{be:1}\cite{bf:1}\cite{bm:1}\cite{km:1}\cite{lt:1}\cite{lt:2}\cite{lt:3}\cite{rt:1}\cite{rt:2}
and have been intensively investigated by several authors
\cite{ehx:1}\cite{ehx:2}\cite{ex:1}\cite{ge:1}\cite{ge:2}\cite{km:2}.
On the other hand, B-twisted sigma models coupled to gravity has not been
understood mathematically (at least for $g>0$ sector) but
topological Landau--Ginzburg models coupled to gravity have been studied based
on singularity theory \cite{lvw:1}\cite{vw:1}\cite{mal:1}.
Recently, Barannikov and Kontsevich have proposed a new construction of
Frobenius manifolds for Calabi--Yau manifolds with arbitrary dimension
\cite{bk:1}, which is conjectured to be a mirror partner for that of
Gromov--Witten theory.
The purpose of this paper is to describe the mathematical foundation of
topological Landau--Ginzburg models coupled to gravity at genus $0$
in terms of the theory of primitive forms which was introduced by K.~Saito
\cite{sa:1}\cite{sa:2} in his study for the period mapping for a universal
unfolding of a function with an isolated critical singularity.
Then, we discuss the mirror symmetry for Calabi--Yau manifolds and $\C\P^1$.
Especially, we can show that the mirror partner of $\C\P^1$ is the
theory of primitive forms associated to a function $f=z+qz^{-1}$ on $\C^*$.
This reproduce the results in \cite{ehx:2}.
In section 1, we review the theory of primitive forms \cite{sa:1}\cite{sa:2}.
Some definitions and notations are modified so that they can be identified
with those in physics literatures.
We will first introduce the notion of a Landau--Ginzburg system $F(z,t)$
over a frame $(Z,X,S,T)$, which corresponds to perturbation theory of the
original Landau--Ginzburg models.
Then we will introduce the chiral ring structure
(the residual product $\circ$),
the perturbed large phase space $\pi_*{\cal H}^{(0)}_F$, the
Gauss--Manin connection $\nabla$ on $\pi_*{\cal H}^{(0)}_F$, the Euler
vector field $E$ and the higher residue pairings $K^{(k)}$.
We will see that a primitive form $\zeta^{(0)}$ associated to $F$ can be
defined if the above data are given.
In section 2, we first review the Frobenius structure defined in \cite{du:1}
\cite{ma:1}.
We will show that the theorem.
\begin{thm*}
Let $F(z,t)$ be an LG system over a frame $(Z,X,S,T)$ and let us assume that a
primitive form $\zeta^{(0)}$ is given.
Then $S$ is a Frobenius manifold with an identity $\delta_0$ and an Euler
vector field $E=w(\delta_0)$.
\end{thm*}
We will construct Cohomological Field Theories \cite{km:1}\cite{ma:1}
associated to primitive forms.
We will see that the large phase space and the gravitational descendants are
given by $H_F=\pi_*{\cal H}^{(0)}_F/m_T\pi_*{\cal H}^{(0)}_F$ and
$\nabla_{\delta_i}\zeta^{(-d-1)}|_{t'=0}$, respectively.
Then we see that
\begin{thm*}
Let $F(z,t)$ be an LG system over a frame $(Z,X,S,T)$.
Then a primitive form $\zeta^{(0)}$ uniquely defines the genus $0$ free energy
$\Phi_F^{grav}$ up to quadratic terms, which is a generating function of
correlation functions of the LG models coupled gravity whose superpotential is
given by $f(z)=F(z,0)$.
\end{thm*}
Note that if we want to have a Frobenius structure, we must fix a primitive
form $\zeta^{(0)}$ for the LG system $F$.
In section 3, we discuss the relations between the primitive forms and the
mirror symmetry.
We will first review the period mapping of the primitive forms.
Then we will show that the canonical coordinates with integral exponents have
a geometrical meaning, i.e., which are one to one correspondence with the
primitive part of the cohomology of $V$.
In particular, if $V$ is a smooth Calabi--Yau manifold,
we can consider that the Frobenius submanifold $S'$ of $S$ (integral exponents
sector) is a Frobenius submanifold of that of Barannikov--Kontsevich
construction from the viewpoint of the LG orbifold theory \cite{iv:1}.
We will also show that a relation between primitive forms and (normalized)
holomorphic forms on Calabi--Yau manifolds.
Finally, we will describe the mirror symmetry of $\C\P^1$ in terms of
primitive forms;
\begin{thm*}
A primitive form for LG system
\begin{equation}
F(z,t)=t^0+z+q\exp(t^1)z^{-1}
\end{equation}
is given by
\begin{equation}
\zeta^{(0)}=\frac{dz}{z}.
\end{equation}
Then we have
\begin{equation}
\Phi_{\C\P^1}^{st}(t)=\Phi_F^{grav}(\tilde{t}),
\end{equation}
where $\Phi_{\C\P^1}^{st}$ is the $($large phase space$)$ stable genus $0$
generating function for Gromov--Witten Invariants
\begin{equation}
\Phi_{\C\P^1}^{st}(t):=\sum_{n\ge 3,(a_i,d_i)}\frac{1}{n!}t^{a_1}_{d_1}\dots
t^{a_n}_{d_n}\sum_{\beta\in \Z_{\ge 0}}q^\beta
\left<\sigma_{d_1}({\cal O}_{a_1})\dots\sigma_{d_n}({\cal O}_{a_n})
\right>_{0,\beta[\C\P^1]}
\end{equation}
and
\begin{equation}
\tilde{t}^i_d:=t^i_d+\sum_{(j,e)}t^j_e\sum_{\beta\in\Z_{\ge 0}}q^\beta
\left<\sigma_{e-d-1}({\cal O}_j){\cal O}^i\right>_{0,\beta[\C\P^1]}.
\end{equation}
This shows that the pair $(\C^*,f)$ is a mirror of $\C\P^1$.
\end{thm*}
This approach may be valid for the other toric Fano manifolds or Grassmannian
manifolds whose quantum cohomology rings are given by the Jacobian ring
\cite{ba:1}\cite{ehx:2}.
We can apply the notion of primitive forms to the Seiberg--Witten theory.
For example, the Seiberg--Witten form $\lambda_{SW}$ for $A_1$ is a
meromorphic $1$-form on
$$
F(x,t^0,\mu)=
W_{A_1}(x,t^0+z+\mu^2 z^{-1})=x^2+t^0+z+\mu^2 z^{-1}=0,
$$
and given by $\lambda_{SW}=x\frac{dz}{z}$.
In this case, the Seiberg--Witten form exactly coincides with the primitive
form ($\zeta^{(-1)}$) for the above defining equation $F(x,t^0,\Lambda)$.
For other gauge groups, we see that $\lambda_{SW}$ is given by
the ``joint'' of the primitive forms for $\C^*$ and $ADE$-type singularity
and reproduces the results of \cite{iy:1}.
We will discuss this in detail in the next paper.
After the completion of this paper, we noticed a preprint by Manin \cite{ma:2}
where the construction of Frobenius manifolds by primitive forms is presented.
\subsection*{Acknowledgement}
~~~
I am deeply grateful to Professor Kyoji Saito for his encouragement.
I also would like to thank Professor Tohru Eguchi, Professor Shinobu Hosono
and Professor Toshiya Kawai for useful discussions.
\tableofcontents
%
\newpage
\section{LG system and primitive form}
\subsection{Landau--Ginzburg system}
~~~
Let us first take the holomorphic function $f$
({\it superpotential}),
$$
f(z):(\C^{n+1},0)\to (\C,0),
$$
which is a quasi--homogeneous function of $z_i$ and has an isolated
singularity at base point $0$.
Then the {\it Jacobian ring} (or {\it chiral ring}) ${\cal R}_0$ defined by
$$
{\cal R}_0:={\cal O}_{\C^{n+1},0}\left/(\frac{\partial f}{\partial z_0},
\dots ,\frac{\partial f}{\partial z_n})\right.
$$
becomes a finite dimensional $\C$-vector space and we denote
$\dim_\C{\cal R}_0$ by $\mu$.
Let us consider a pair $(S,\delta_0 ,0)$ where $S$ is a $\mu$-dimensional
complex manifold with a base point $0\in S$ and $\delta_0$
is a holomorphic vector field on $S$ which is non-singular at $0$.
There exist a $(\mu-1)$-dimensional manifold $T$ with a base point $0\in T$
and a holomorphic submersion $\pi$ defined locally,
\begin{equation}
\pi :(S,0)\to (T,0)\ such~that~\ \pi^{-1}{\cal O}_T=\{ g\in{\cal O}_S:\delta_0
g=0\}
\end{equation}
where ${\cal O}_S$ and ${\cal O}_T$ are structure sheaves of $S$ and $T$.
By choosing neighborhoods of base points of $S$ and $T$,
we may assume that $T$ is the orbit space of $\delta_0$.
Let us denote by $Der_S$ the sheaf of
germs of holomorphic vector fields on $S$ .
Put
\begin{equation}
{\cal G} := \{ \delta\in \pi_* Der_S :[\delta_0,\delta ]=0\}
\end{equation}
where $[\ ,\ ]$ is the bracket product of vector fields.
The module ${\cal G}$ is naturally ${\cal O}_T$-free module of rank $\mu$
closed under the Lie-bracket.
\begin{rem}
We have the following exact sequence of ${\cal O}_T$-Lie algebras:
\begin{equation}
0\to{\cal O}_T\delta_0\to{\cal G}\to Der_T\to 0.
\end{equation}
In this situation, the following are equivalent.
\begin{enumerate}
\item To give a holomorphic function $t^0$ on $(S,\delta_0,0)$ such that
$t^0(0)=0$ and $\delta_0 t^0=1$.
\item To give a splitting of the above exact sequence as an
${\cal O}_T$-Lie algebra.
\end{enumerate}
\end{rem}
Let $X$ be an $(n+\mu)$-dimensional complex manifold with a base point
$0\in X$, and let $q:(X,0)\to (T,0)$ be a smooth surjective map.
Let $(Z,0)$ be the fiber product with the diagram,
\begin{equation}
\begin{CD}
(Z,0) @>\hat{\pi}>> (X,0) \\
@V p VV @VV q V \\
(S,0) @>\pi>> (T,0)
\end{CD}
\end{equation}
There exists a unique lifting $\hat{\delta_0}$ on $Z$ of the primitive vector
field $\delta_0$ on $S$ such that $\hat{\pi}:(Z,0)\to (X,0)$ is the
projection to the orbit space of $\hat{\delta_0}$.
We shall call the above diagram a {\bf frame},
which will often be denoted simply by $(Z,X,S,T)$.
\begin{rem}
The following are equivalent.
\begin{enumerate}
\item To give a holomorphic function $F(z,t)$ on $Z$
such that $\hat{\delta_0}F=1$ and $F(0)=0$.
\item To give a holomorphic mapping $\varphi :(X,0)\to (S,0)$,
such that $\pi\circ\varphi=q$.
\end{enumerate}
\end{rem}
\begin{defn}
A holomorphic function $F(z,t)$ on $Z$ is called a
{\it Landau--Ginzburg system} ({\it LG system}) over a frame $(Z,X,S,T)$
if it satisfies the following three conditions:
\begin{enumerate}
\item $F(z,0)=f(z)$.
\item $\hat{\delta_0}F=1$, $F(0)=0$.
\item The correspondence
\begin{equation}
T(S)_0\to {\cal R}_0={\cal O}_{p^{-1}(0),0}\left/
(\frac{\partial f}{\partial z_0},\dots ,\frac{\partial f}{\partial z_n})
\right.,~~~
\delta\mapsto\hat{\delta}F|_{p^{-1}(0)},
\end{equation}
is bijective.
\end{enumerate}
\end{defn}
%
\subsection{The residual product $\circ$}
~~~
Let $F(z,t)$ be an LG system over a frame $(Z,X,S,T)$ and let
$\varphi :(X,0)\to (S,0)$ be an associated holomorphic map.
Let ${\cal C}$ be the subvariety of the critical points
of $\varphi$ in $X$ defined by the ideal
$(\frac{\partial F}{\partial z_0},\dots ,
\frac{\partial F}{\partial z_n}){\cal O}_X$.
Then we see that the restriction $q|_{\cal C}:{\cal C}\to T$ is a
$\mu$-sheeted branched covering and $q_*{\cal O}_{\cal C}$ is an
${\cal O}_T$-free module of rank $\mu$.
Let us put $\hat{{\cal C}}:=\hat{\pi}^{-1}{\cal C}\subset Z$,
which is defined by the ideal $(\frac{\partial F}{\partial z_0},\dots ,
\frac{\partial F}{\partial z_n}){\cal O}_Z$.
Viewing $X$ as a hypersurface in $Z$ defined by $F=0$,
we have the following exact sequence:
\begin{equation}
0\to p_* {\cal O}_{\hat{{\cal C}}}\stackrel {p_* F}{\to } p_*
{\cal O}_{\hat{{\cal C}}}\to \varphi_* {\cal O}_{\cal C}\to 0,
\end{equation}
where $p_* F$ means an operator of a multiplication of $F$.
\begin{lem}
Let $F(z,t)$ be an LG system over a frame $(Z,X,S,T)$.
Then the ${\cal O}_S$-homomorphism
\begin{equation}
Der_S\to p_*{\cal O}_{\hat{{\cal C}}},~~~
\delta\mapsto \hat{\delta }F|_{\hat{{\cal C}}},
\end{equation}
and the ${\cal O}_T$-homomorphism
\begin{equation}
{\cal G}\to q_*{\cal O}_{\cal C},~~~
\delta\mapsto \hat{\delta }F|_{\cal C}.
\end{equation}
are bijections.
\end{lem}
\begin{cor}
We can give a ring structure and a $t^0$-multiplication structure on ${\cal G}$
by,
\begin{equation}
\hat{(\delta\circ\delta')}F|_{\cal C}:=
\hat{\delta}F|_{\cal C}\cdot\hat{\delta'}F|_{\cal C}
\end{equation}
and
\begin{equation}
(\hat{t^0\circ\delta})F|_{\cal C}:=t^0|_{\cal C}\cdot\hat{\delta}F|_{\cal C}.
\end{equation}
\end{cor}
\begin{defn}
We call the above $\circ$ the {\it residual product}.
\end{defn}
Since $p|_{\hat{{\cal C}}}:\hat{{\cal C}}\to S$ is $\mu$-sheeted
branched covering, $p_* {\cal O}_{\hat{{\cal C}}}$ is an ${\cal O}_S$-free
module of rank $\mu$, and hence the sequence is an ${\cal O}_S$-free
resolution of $\varphi_* {\cal O}_{\cal C}$.
The determinant of the ${\cal O}_S$-endomorphism $p_* F$
\begin{equation}
\Delta :=\det (p_* F)
\end{equation}
is called the {\it discriminant}.
${\cal D}:=\varphi ({\cal C})$ is defined by
the ideal sheaf ${\cal J}=(\Delta )$.
Let us put
\begin{equation}
Der_S(-\log {\cal D}):=\{ \delta\in Der_S~|~\delta(\Delta)\in(\Delta)\}.
\end{equation}
which is an ${\cal O}_S$-free module of rank $\mu$, closed under bracket
product. Considering the ${\cal O}_T$-homomorphism
\begin{equation}
r:\pi_*Der_S\to q_*{\cal O}_{\cal C},
\end{equation}
we see that $\ker(r)=\pi_*Der_S(-\log {\cal D})$ and
we have a direct sum decomposition
\begin{equation}
\pi_*Der_S={\cal G}\oplus\pi_*Der_S(-\log {\cal D}),
\end{equation}
and an ${\cal O}_T$-homomorphism $w$
\begin{equation}
w:{\cal G}\to \pi_*Der_S(-\log {\cal D}),~~~
\delta\mapsto w(\delta)=t^0\delta-t^0\circ\delta.
\end{equation}
Let us define the specific vector field $E:=w(\delta_0)$ which we call the
{\it Euler vector field}.
The discriminant $\Delta$ is homogeneous degree $\mu$, with respect to E, i.e.,
\begin{equation}
E\Delta=\mu\Delta.
\end{equation}
%
\subsection{De Rham cohomology ${\cal H}^{(0)}_F$}
~~~
Let $\Omega^p_X$ be the sheaf of germs of holomorphic $p$-forms on $X$,
$\Omega^p_{X/T}:=\Omega^p_X/dq^*\Omega^1_T\wedge\Omega^{p-1}_X$
and $\Omega^p_{X/S}:=\Omega^p_X/d\varphi^*\Omega^1_S\wedge\Omega^{p-1}_X$
be the sheaves of germs of holomorphic relative p-forms with respect to $q$
or $\varphi$ respectively.
In particular let us define,
\begin{equation}
\Omega_F:=\Omega^{n+1}_{X/S}.
\end{equation}
Note that if we fix the holomorphic $n+1$-form $dz_0\wedge\dots\wedge dz_n$,
we have $\Omega_F\simeq{\cal O}_{\cal C}dz_0\wedge\dots\wedge dz_n$.
Let us define the ${\cal O}_S$-free modules,
\begin{equation}
\begin{aligned}
{\cal H}^{(0)}_F :&=\varphi_*\Omega^{n+1}_{X/T}/dF^0\wedge
d(\varphi_*\Omega^{n-1}_{X/T})\\
{\cal H}^{(-1)}_F :&=\varphi_*\Omega^{n}_{X/T}/dF^0\wedge
(\varphi_*\Omega^{n-1}_{X/T})+d(\varphi_*\Omega^{n-1}_{X/T})\\
&=\varphi_*\Omega^{n}_{X/S}/d(\varphi_*\Omega^{n-1}_{X/S})\\
{\cal H}^{(-2)}_F :&=
\ker(d:\varphi_*\Omega^{n}_{X/S}\to\varphi_*\Omega^{n+1}_{X/S})
/d(\varphi_*\Omega^{n-1}_{X/S})\\
&=\R^{n}\varphi_*(\Omega^{\cdot}_{X/S},d),
\end{aligned}
\end{equation}
where $F^0$ is a function on $X$ defined by $F=t^0-F^0$.
They are ${\cal O}_S$-free modules of rank $\mu$ because of
the Milnor's fibration theorem:
\begin{thm}[Milnor's fibration theorem]
The restriction
\begin{equation}
\varphi|_{X-\varphi^{-1}({\cal D})}:X-\varphi^{-1}({\cal D})\to S-{\cal D}
\end{equation}
is a locally trivial fibration, whose general fiber $X_s=\varphi^{-1}(s)$ for
$s\in S-{\cal D}$ is homeomorphic to a bouquet of $n$-spheres.
\end{thm}
Indeed, we see that
\begin{equation}
H^i(X_s,\C)\simeq
\begin{cases}
\C,& i=0,\\
0,& i\ne 0~{\rm and}~n,\\
\C^\mu,& i=n.
\end{cases}
\end{equation}
The exterior differentiation $d:\Omega^n_{X/T}\to\Omega^{n+1}_{X/T}$
induces an ${\cal O}_S$-exact sequence,
\begin{equation}
0\to {\cal H}^{(-2)}_F\hookrightarrow {\cal H}^{(-1)}_F\stackrel{d}{\to }
\varphi_*\Omega_F\to 0.
\end{equation}
The exterior product $dF^0:\Omega^n_{X/T}\to \Omega^{n+1}_{X/T}$
induces an ${\cal O}_S$-exact sequence,
\begin{equation}
0\to {\cal H}^{(-1)}_F\stackrel{dF^0}{\hookrightarrow }{\cal H}^{(0)}_F
\to \varphi_*\Omega_F\to 0.
\end{equation}
\begin{defn}
The covariant differentiation
\begin{equation}
\nabla :Der_S\times {\cal H}^{(-1)}_F\to {\cal H}^{(0)}_F,
\end{equation}
defined by
\begin{equation}
\nabla_{\frac{\partial }{\partial t^i}}[\zeta]=(-1)^i[(dt^1\wedge\dots
\wedge dt^{\mu-1})^{-1}dt^0\wedge\dots\wedge\hat{dt^i}\wedge\dots\wedge dt^{\mu-1}
\wedge d\zeta],
\end{equation}
is called the {\it Gauss--Manin connection}.
\end{defn}
\begin{prop}
The Gauss--Manin connection is integrable and
\begin{equation}
\nabla_{\delta_0}:{\cal H}^{(-1)}_F\simeq {\cal H}^{(0)}_F,~~~
\nabla_{\delta_0}[\zeta]=[d\zeta].
\end{equation}
Then the map
\begin{equation}
\nabla^{-1}_{\delta_0}:{\cal H}^{(0)}_F\to{\cal H}^{(0)}_F,~~~
[\zeta]\mapsto[dF^0\wedge d^{-1}\zeta],
\end{equation}
is well--defined.
\end{prop}
Thus we can consider the following decreasing sequence
\begin{equation}
{\cal H}^{(-k-1)}_F:=\{ \omega\in {\cal H}^{(-k)}_F:\nabla_{\delta_0}
\omega\in {\cal H}^{(-k)}_F\},~~~k\in \N
\end{equation}
of ${\cal O}_S$-free submodules of ${\cal H}^{(0)}_F$.
We get the ${\cal O}_S$-exact sequences,
\begin{equation}
0\to {\cal H}^{(-k-1)}_F\hookrightarrow {\cal H}^{(-k)}_F
\stackrel{r^{(k)}}{\to }\varphi_*\Omega_F\to 0,~~~k\in\N,
\end{equation}
and the Gauss-Manin connection $\nabla$,
\begin{equation}
\nabla :Der_S\times{\cal H}^{(-k-1)}_F\to{\cal H}^{(-k)}_F,~~~k\in \N.
\end{equation}
\begin{rem}
Since $r^{(k)}(\nabla_\delta\zeta)=r(\delta)r^{(k)}(\nabla_{\delta_0}\zeta)$,
we have
\begin{equation}
\nabla:Der_S(-\log({\cal D}))\times{\cal H}^{(-k)}_F\to{\cal H}^{(-k)}_F,
~~~k\in\N.
\end{equation}
This means that the Gauss--Manin connection is regular singular.
\end{rem}
%
\subsection{Higher residue pairings $K^{(k)}$ and
primitive form $\zeta^{(0)}$}
~~~
We will consider the direct images $\pi_*{\cal H}^{(-k)}_F$.
Note that they are not finite ${\cal O}_T$-modules.
The pairing $K^{(k)}$ on $\pi_*{\cal H}^{(0)}_F$ in the following theorem
is called the {\it higher residue pairings} \cite{sa:1}\cite{sa:2}.
\begin{thm}
Let $F(z,t)$ be an LG system over a frame $(Z,X,S,T)$.
Then there exists an infinite sequence of ${\cal O}_T$-bilinear forms,
\begin{equation}
K^{(k)}:\pi_*{\cal H}^{(0)}_F\times\pi_*{\cal H}^{(0)}_F\to {\cal O}_T,~~~
k\in\N
\end{equation}
with the following properties$:$
\begin{enumerate}
\item $K^{(k)}$ is symmetric $($skew-symmetric$)$ when $k$ is even $($odd $)$,
respectively.
\item
$K^{(0)}([\phi_1dz],[\phi_2dz])={\rm Res}_{X/T}
\left[
\begin{gathered}
\phi_1\phi_2dz_0\wedge\dots\wedge dz_n\\
\frac{\partial F}{\partial z_0}\dots\frac{\partial F}{\partial z_n}
\end{gathered}
\right]$, \\
for $[\phi_1dz_0\wedge\dots\wedge dz_n]$,
$[\phi_2dz_0\wedge\dots\wedge dz_n]\in\pi_*{\cal H}^{(0)}_F$.
\item
$K^{(k)}(\omega_1,\omega_2)=K^{(k-1)}(\nabla_{\delta_0}\omega_1,\omega_2)$,
for $\omega_1\in\pi_*{\cal H}^{(-1)}_F$, $\omega_2\in\pi_*{\cal H}^{(0)}_F$.
\item
$K^{(k)}(t^0\omega_1,\omega_2)-K^{(k)}(\omega_1,t^0\omega_2)=(n+k)K^{(k-1)}
(\omega_1,\omega_2)$, \\
for $\omega_1$, $\omega_2\in\pi_*{\cal H}^{(0)}_F$.
\item
$\delta K^{(k)}(\omega_1,\omega_2)=K^{(k)}(\nabla_\delta\omega_1,\omega_2)+
K^{(k)}(\omega_1,\nabla_\delta\omega_2)$, \\
for $\omega_1$, $\omega_2\in\pi_*{\cal H}^{(-1)}_F,~~\delta\in{\cal G}$.
\end{enumerate}
Such a $K^{(k)}$ with the above properties is unique up to constant factor.
\end{thm}
\begin{rem}
$K^{(0)}$ induces a non-degenerate ${\cal O}_T$-bilinear form $J$:
\begin{equation}
J:q_*\Omega_F\times q_*\Omega_F\to {\cal O}_T.
\end{equation}
\end{rem}
\begin{defn}\label{de:pr}
An element $\zeta^{(0)}\in\Gamma(S,{\cal H}^{(0)}_F)$
is called a {\it primitive form} if it satisfies the following conditions:
\begin{enumerate}
\item Invertibility. $r^{(0)}(\zeta^{(0)})\in\Gamma(C,\Omega_F)$ is an
${\cal O}_{\cal C}$-free base of $\Omega_F$.
\item First integrability. $K^{(1)}(\nabla_\delta\zeta^{(-1)},
\nabla_{\delta'}\zeta^{(-1)})=0$ for $\delta,\delta'\in{\cal G}$.
\item Homogeneity. $\nabla_E\zeta^{(0)}=(r-1)\zeta^{(0)}$ for a constant $r$.
\item $K^{(k)}(\nabla_\delta\nabla_{\delta'}\zeta^{(-2)},\nabla_{\delta''}
\zeta^{(-1)})=0$ for $k\ge 2$ and $\delta,\delta',\delta''\in{\cal G}$.
\item $K^{(k)}(t^0\nabla_\delta\zeta^{(-1)},\nabla_{\delta'}\zeta^{(-1)})=0$
for $k\ge 2$ and $\delta,\delta'\in{\cal G}$.
\end{enumerate}
where $\zeta^{(-k)}:=(\nabla_{\delta_0})^{-k}\zeta$.
\end{defn}
\begin{thm}
If $S$ is small enough, then there exists a primitive form $\zeta^{(0)}$.
\end{thm}
%
\subsection{Canonical coordinate system}
Since $r^{(0)}(\zeta^{(0)})$ is invertible, there is an
${\cal O}_T$-isomorphism
\begin{equation}
{\cal G}\simeq q_*\Omega_F,~~~\delta\mapsto\delta F|_{\cal C}\cdot
r^{(0)}(\zeta^{(0)}).
\end{equation}
Then we can define a connection
\begin{equation}
\tilde{\nabla}:Der_T\times{\cal G}\to {\cal G}
\end{equation}
by the relation
\begin{equation}
\nabla_\delta\nabla_{\delta'}\zeta^{(-2)}=\nabla_{\delta\circ\delta'}
\zeta^{(-1)}+\nabla_{\tilde{\nabla}_\delta\delta'}\zeta^{(-2)},
\end{equation}
which depends only on $r^{(0)}(\zeta^{(0)})$.
The integrability of the Gauss--Manin connection $\nabla$ on
${\cal H}^{(0)}_F$ follows that
\begin{equation}\label{eq:nab}
\tilde{\nabla}_\delta(\delta'\circ\delta'')-\delta'\circ\tilde{\nabla}_\delta
\delta''-\tilde{\nabla}_{\delta'}(\delta\circ\delta'')+\delta\circ
\tilde{\nabla}_{\delta'}\delta''=[\delta,\delta']\circ\delta'',
~~\delta,\delta',\delta''\in{\cal G},
\end{equation}
and
\begin{equation}
[\tilde{\nabla}_\delta,\tilde{\nabla}_{\delta'}]=
\tilde{\nabla}_{[\delta,\delta']},~~~\delta,\delta'
\in{\cal G}.
\end{equation}
By putting $\delta''=\delta_0$ in \eqref{eq:nab},
we see that this connection is torsion free
\begin{equation}
\tilde{\nabla}_\delta\delta'-\tilde{\nabla}_{\delta'}\delta=[\delta,\delta'],
~~~\delta,\delta'
\in{\cal G}.
\end{equation}
From the fifth property of $K^{(k)}$, we see that this connection is metric
with respect to $J$,
\begin{equation}
\delta J(\delta',\delta'')=J(\tilde{\nabla}_\delta\delta',\delta'')
+J(\delta',\tilde{\nabla}_\delta\delta''),
~~~\delta,~\delta',~\delta''\in{\cal G},
\end{equation}
where we denote $J(\delta'F|_{\cal C}\cdot r^{(0)}(\zeta^{(0)}),
\delta''F|_{\cal C}\cdot r^{(0)}(\zeta^{(0)}))$ by $J(\delta,\delta')$
for short.
We can define an ${\cal O}_T$-endomorphism $N$
\begin{equation}
N:{\cal G}\to{\cal G}
\end{equation}
by the relation
\begin{equation}
t^0\nabla_\delta\zeta^{(-1)}=\nabla_{t^0\circ\delta}\zeta^{(-1)}
+\nabla_{N\delta}\zeta^{(-2)},~~~\delta\in{\cal G}.
\end{equation}
The commutation relation $[\delta,t^0]=\delta t^0$ follows
\begin{equation}
\tilde{\nabla}_\delta(t^0\circ\delta')+\delta\circ(N+1)\delta'=N(\delta\circ
\delta')+t^0\circ\tilde{\nabla}_\delta\delta'+(\delta t^0)\delta',
~~~\delta,\delta'\in{\cal G},
\end{equation}
and
\begin{equation}\label{eq:N}
\tilde{\nabla}_\delta(N\delta')=N(\tilde{\nabla}_\delta\delta'),
~~~\delta,\delta'\in{\cal G}.
\end{equation}
From the fourth property of $K^{(k)}$, we get
\begin{equation}
J(N\delta,\delta')+J(\delta,N\delta')=(n+1)J(\delta,\delta'),
~~~\delta,\delta'\in{\cal G},
\end{equation}
i.e.
\begin{equation}\label{eq:dual}
N+N^*=(n+1)id,
\end{equation}
where $N^*$ is the adjoint of $N$ with respect to the metric $J$.
Since the connection $\tilde{\nabla}$ on ${\cal G}$ is torsion free and
integrable, there exists a unique coordinate system (called a {\it flat
coordinate system}) $t^0,\dots ,t^{\mu-1}$ of $S$ up to linear transformations
such that the horizontal vector space of $\tilde{\nabla}$ is given by
\begin{equation}
\ker\tilde{\nabla} =\bigoplus^{\mu-1}_{i=0}\C\frac{\partial }{\partial t^i}.
\end{equation}
The equation \eqref{eq:N} implies that $N$ is a $\C$-linear endomorphism of
$\ker\tilde{\nabla}$.
\begin{defn}
The eigenvalues $\{\alpha_0,\alpha_1,\dots ,\alpha_{\mu-1}\}$ of
$N$ is called the {\it exponents}.
\end{defn}
We introduce the following {\it Poincar\'e polynomial} associated to the
LG system $F$:
\begin{equation}
\chi(F)(y,\bar{y}):=\sum_{i=0}^{\mu-1}(y\bar{y})^{q_{i}},
\end{equation}
where $q_i:=\alpha_i-r$, $i=0,\dots,\mu-1$.
Since we can choose exponents such that
$\alpha_0(:=r)\le\dots\le\alpha_{\mu-1}$
and exponents have the following duality property:
\begin{equation}
\{\alpha_0,\dots ,\alpha_{\mu-1}\}=
\{n+1-\alpha_0,\dots ,n+1-\alpha_{\mu-1}\},
\end{equation}
it is obvious from \eqref{eq:dual} that the duality property of the
Poincar\'e polynomial,
\begin{equation}
\chi(F)(y,\bar{y})=(y\bar{y})^{\hat{c}_F}\chi(F)(y^{-1},\bar{y}^{-1}).
\end{equation}
where $\hat{c}_F:=q_{\mu-1}=(n+1)-2r$.
Let $t^0,\dots ,t^{\mu-1}$ be a flat coordinate system.
If $N\left(\frac{\partial}{\partial t^i}\right)=
\alpha_i\frac{\partial}{\partial t^i}$, then we can define
the {\it degree} of $t^i$ by
\begin{equation}
\deg(t^i):=1-q_i,~~~i=0,\dots,\mu-1.
\end{equation}
\begin{defn}
A coordinate system $t^0,\dots ,t^{\mu-1}$ of $S$ is called a
{\it canonical coordinate system} if it is a flat coordinate system and
\begin{equation}
N\left(\frac{\partial}{\partial t^i}\right)=
\alpha_i\frac{\partial}{\partial t^i},~~~ i=1,\dots ,\mu-1.
\end{equation}
\end{defn}
\begin{rem}
In a canonical coordinate system, the Euler vector field $E=w(\delta_0)$
is given by
\begin{equation}
E=\sum_{i=0}^{\mu-1}\deg(t^i)t^i\frac{\partial}{\partial t^i}.
\end{equation}
\end{rem}
\begin{rem}
This canonical coordinates coincides with the one introduced by \cite{bcov:1}.
\end{rem}
\section{Cohomological Field Theory associated to LG system}
\subsection{Frobenius structure}
~~~
In this subsection, we will show that $S$ is a Frobenius manifold.
We first review the Frobenius structure \cite{du:1}\cite{ma:1}.
\begin{defn}
An affine flat structure on the supermanifold $M$ is a subsheaf
${\cal T}^f_M\subset{\cal T}_M$ of linear spaces of pairwise (super)commuting
vector fields, such that ${\cal T}_M={\cal T}^f_M\otimes {\cal O}_M$
(tensor product over the ground field).
Sections of ${\cal T}^f_M$ is called flat vector field.
\end{defn}
A Riemannian metric on $M$ is an even symmetric pairing
$\eta:S^2({\cal T}_M)\to{\cal O}_M$, inducing an isomorphism
$\eta':{\cal T}_M\to {\cal T}^*_M$.
We put $\eta_{ab}:=\eta(\partial_a,\partial_b)$.
\begin{defn}
A Riemannian metric $\eta\in{\cal T}^*_M\otimes{\cal T}^*_M$ is compatible with
the structure ${\cal T}^f_M$
if $\eta(X,Y)$ is constant for all $X,Y\in{\cal T}^f_M$.
\end{defn}
\begin{defn}
Let $M$ be a supermanifold. Consider a triple $({\cal T}^f_M,\eta,C)$
consisting of an affine flat structure, a compatible metric and an even
symmetric tensor $C:S^3({\cal T}_M)\to {\cal O}_M$.
Define an ${\cal O}_M$-bilinear symmetric multiplication $\circ$ on
${\cal T}_M$:
\begin{equation}
\begin{CD}
{\cal T}_M\otimes{\cal T}_M\; @>>>\; S^2({\cal T}_M)\; @>{C'}>>\; {\cal T}^*_M
\; @>\eta'>>\; {\cal T}_M,
\end{CD}
\end{equation}
where primes denotes a partial dualization.
Then $M$ endowed with this structure is called a {\it pre--Frobenius} manifold.
\end{defn}
A local even function $\Phi\in{\cal O}_M$ is called a
{\it Frobenius potential} if
\begin{equation}
C(X,Y,Z)=(XYZ)\Phi,~~~X,Y,Z\in{\cal T}_M.
\end{equation}
\begin{defn}
A pre--Frobenius manifold is called {\it Frobenius}, if $C$ everywhere
locally admits a potential and the multiplication $\circ$ is associative.
\end{defn}
\begin{rem}
If a Frobenius potential exists, it is unique up to a quadratic polynomial in
flat local coordinates.
\end{rem}
\begin{rem}
In flat local coordinates $t^a$, then
\begin{equation}
C_{abc}=\partial_a\partial_b\partial_c\Phi,
\end{equation}
\begin{equation}
\partial_a\circ\partial_b=\sum_cC_{ab}^c\partial_c,
\end{equation}
where
\begin{equation}
C_{ab}^c:=\sum_dC_{abd}\eta^{cd},~~~~(\eta^{ab}):=(\eta_{ab})^{-1}.
\end{equation}
\end{rem}
\begin{defn}
Let $(M,\eta,\Phi)$ be a Frobenius manifold.
An even vector field $e\in{\cal T}_M$ is called a {\it identity} if
$e\circ X=X$ for all vector fields $X\in{\cal T}_M$.
\end{defn}
\begin{defn}
An even vector field $E\in{\cal T}_M$ is called an {\it Euler vector field}
if it satisfies the following:
\begin{equation}
E(\eta(X,Y))-\eta([E,X],Y)-\eta(X,[E,Y])=D\eta(X,Y),~~~X,Y\in{\cal T}_M,
\end{equation}
and
\begin{equation}
[E,X\circ Y]-[E,X]\circ Y-X\circ [E,Y]=d_0X\circ Y,~~~X,Y\in{\cal T}_M.
\end{equation}
\end{defn}
\begin{thm}
Let $F(z,t)$ be an LG system over a frame $(Z,X,S,T)$ and let us assume that a
primitive form $\zeta^{(0)}$ is given.
Then $S$ is a Frobenius manifold with an identity $\delta_0$ and an Euler
vector field $E=w(\delta_0)$.
\end{thm}
\begin{pf}
It is clear that a canonical coordinate system $t^i$, $i=1,\dots ,\mu$ give an
affine flat structure. We will show the existence of a compatible metric
$\eta$ and a symmetric tensor $C$.
Let
\begin{equation}
\eta_F(\delta,\delta'):=\sum_{i,j=0}^{\mu-1}f^ig^jJ(\delta_i,\delta_j),
\end{equation}
for $\delta=\sum_{i=0}^{\mu-1}f^i\delta_i$,
$\delta'=\sum_{i=0}^{\mu-1}f^j\delta_j\in Der_S$.
Then $\eta_F(\delta_i,\delta_j)$ is constant for all flat vector fields
$\delta_i,\delta_j$ because of the second property of the definition
\ref{de:pr};
\begin{equation}
\begin{aligned}
\delta_k\eta_F(\delta_i,\delta_j) &= \delta_k
K^{(0)}(\nabla_{\delta_i}\zeta^{(-1)},\nabla_{\delta_j}\zeta^{(-1)})\\
&= K^{(1)}(\nabla_{\delta_k}\nabla_{\delta_i}\zeta^{(-2)},
\nabla_{\delta_j}\zeta^{(-1)})-K^{(1)}(\nabla_{\delta_i}\zeta^{(-1)},
\nabla_{\delta_k}\nabla_{\delta_j}\zeta^{(-2)})\\
&= K^{(1)}(\nabla_{\delta_k\circ\delta_i}\zeta^{(-1)},
\nabla_{\delta_j}\zeta^{(-1)})-K^{(1)}(\nabla_{\delta_i}\zeta^{(-1)},
\nabla_{\delta_k\circ\delta_j}\zeta^{(-1)})\\
&= 0
\end{aligned}
\end{equation}
Let us define $C$ by the residual product $\circ$
\begin{equation}\label{eq:aso}
\delta_i\circ\delta_j=\sum_{k=0}^\mu C_{ij}^k\delta_k.
\end{equation}
Since
\begin{equation}
\begin{split}
\partial_i C_{jkl}&=\delta_i K^{(0)}(\nabla_{\delta_j}\nabla_{\delta_k}
\zeta^{(-2)},\nabla_{\delta_l}\zeta^{(-1)})\\
&=K^{(0)}(\nabla_{\delta_i}\nabla_{\delta_j}\nabla_{\delta_k}\zeta^{(-2)},
\nabla_{\delta_l}\zeta^{(-1)})\\
&=K^{(0)}(\nabla_{\delta_j}\nabla_{\delta_i}\nabla_{\delta_k}\zeta^{(-2)},
\nabla_{\delta_l}\zeta^{(-1)})\\
&=\delta_j K^{(0)}(\nabla_{\delta_i}\nabla_{\delta_k}
\zeta^{(-2)},\nabla_{\delta_l}\zeta^{(-1)})\\
&=\partial_j C_{ikl},
\end{split}
\end{equation}
there exists a potential $\Phi_F$ such that $C_{ijk}=\partial_i\partial_j
\partial_k\Phi_F$. The associativity is obvious from the definition of the
residual product $\circ$. The rest is clear.
\end{pf}\qed
%
\subsection{Cohomological Field Theory}
~~~
Let $k$ be a supercommutative $\Q$-algebra,
$H$ a ($\Z_2$-graded) free $k$-module of finite rank, $\eta:H\times H\to k$
an even symmetric non--degenerate pairing on $H$.
We introduce formally the {\it large phase space} as linear infinite
dimensional formal supermanifold $\oplus_{d\ge 0}H[d]$ with basis
$\sigma_d({\cal O}_a)$ and coordinates $t^a_d$.
Put $t^a:=t^a_0$, $t:=\{t^a\}$.
A {\it Frobenius potential} on $(H,\eta)$ is given by
an even potential $\Phi\in k[[t]]$, defined up to quadratic terms,
whose third derivatives form the structure constants of the
supercommutative, associative $k[[t]]$-algebra $H\otimes k[[t]]$.
\begin{defn}
A triple $M=(H,\eta,\Phi)$ is called a {\it formal Frobenius manifold}
(over $k$).
\end{defn}
Let $\overline{\cal M}_{0,n}$ be the Deligne--Mumford compactification of the
moduli space of genus $0$ Riemann surfaces with $n$-marked points.
The pair $k$-vector space $H$ and maps $I^M_n$ in the following proposition is
called the (tree level) {\it Cohomological Field Theory}
\cite{km:1}\cite{ma:1}.
\begin{prop}{\rm \cite{km:2}}
For any formal Frobenius manifold $M$, there exists a unique sequence
of $k$-linear maps
\begin{equation}
I^M_n:H^n\to H^*(\overline{\cal M}_{0,n},k),~~~n\ge 3,
\label{eq:co}
\end{equation}
satisfying the following properties$:$
\begin{enumerate}
\item ${\cal S}_n$-covariance $($with respect to the natural action of
${\cal S}_n$ on both sides of \eqref{eq:co}$)$.
\item Splitting, or compatibility with restriction to the boundary divisors$:$
for any stable ordered partition $\sigma:\{1,\dots ,n\}=S_1\coprod S_2$,
$n_i=|S_i|$, and the respective map
\begin{equation}
\varphi_\sigma:\overline{\cal M}_{0,n_1+1}\times\overline{\cal M}_{0,n_2+1}\to
\overline{\cal D}(\sigma)\subset\overline{\cal M}_{0,n}
\end{equation}
we have
\begin{multline}
\varphi^*_\sigma(I^M_n({\cal O}_{a_1}\otimes\dots\otimes{\cal O}_{a_n}))\\
=\epsilon(\sigma)(I^M_{n_1+1}\otimes I^M_{n_2+1})
((\bigotimes_{p\in S_1}{\cal O}_{a_p})\otimes\Delta\otimes
(\bigotimes_{q\in S_2}{\cal O}_{a_q}))
\end{multline}
where $\Delta:=\sum_{a,b}{\cal O}_{a}\eta^{ab}{\cal O}_{b}$,
and $\epsilon(\sigma)$ is the sign of the permutation
induced on the odd arguments ${\cal O}_{a_1},\dots ,{\cal O}_{a_n}$.
\item
\begin{equation}
\int_{\overline{\cal M}_{0,n}}I^M_n({\cal O}_{a_1}\otimes\dots\otimes
{\cal O}_{a_n})=\partial_{a_1}\dots\partial_{a_n}\Phi|_{t=0},
\end{equation}
where the integral denotes the value of the top degree term of $I^M_n$ on the
fundamental class $\overline{\cal M}_{0,n}$.
\end{enumerate}
\end{prop}
Now we can define the correlators with gravitational descendants
$\sigma_d({\cal O}_a)$ by
\begin{multline}
\left<\sigma_{d_1}({\cal O}_{a_1})\dots\sigma_{d_n}({\cal O}_{a_n})\right>\\
:=\int_{\overline{\cal M}_{0,n}}I^M_n({\cal O}_{a_1}\otimes\dots
\otimes{\cal O}_{a_n})c_1({\cal L}_1)^{d_1}\dots c_1({\cal L}_n)^{d_n},
\end{multline}
where ${\cal L}_i$ is the line bundle which has the geometric fiber
$T^*_{x_i}C$ at the point $(C;x_1,\dots ,x_n)$.
\begin{defn}
A formal function $\Phi^{grav}$ on large phase space $\oplus_{d\ge 0}H[d]$
defined by
\begin{equation}
\Phi^{grav}:=\sum_{n,(d_1,a_1),\dots ,(d_n,a_n)}\epsilon(a)\frac{1}{n!}
t^{a_1}_{d_1}\dots t^{a_n}_{d_n}\left<\sigma_{d_1}
({\cal O}_{a_1})\dots\sigma_{d_n}({\cal O}_{a_n})\right>
\end{equation}
is called a {\it genus $0$ free energy}.
\end{defn}
Let us introduce $\mu$-dimensional vector spaces,
\begin{align}
H_F^{(d)} & := \bigoplus_{i=0}^{\mu-1}\C
\nabla_{\delta_i}\zeta^{(-d-1)}\left|_{t'=0}\right. \\
& \simeq
\left(\pi_*{\cal H}^{(-d)}_F\left/\pi_*{\cal H}^{(-d-1)}_F\right.
\right)\left/m_0\left(\pi_*{\cal H}^{(-d)}_F\left/\pi_*{\cal H}^{(-d-1)}_F
\right.\right)\right. ,~~~d\in\N,
\end{align}
where $m_0=(t'):=(t^1,\dots ,t^{\mu-1})$ is the maximal ideal of
${\cal O}_{T,0}$.
\begin{lem}
Let $F(z,t)$ be an LG system over a frame $(Z,X,S,T)$ and let us assume that a
primitive form $\zeta^{(0)}$ is given.
Then $M_F=(H_F^{(0)},\eta_F,\Phi_F)$ is a formal Frobenius manifold over $\C$.
\end{lem}
\begin{pf}
Since $S$ is a Frobenius manifold, it is clear that $M_F$ is a formal Frobenius
manifold over $\C$.
\end{pf}\qed
\begin{rem}
The large phase space is given by
\begin{equation}
H_F:=\pi_*{\cal H}^{(0)}_F\left/m_0\pi_*{\cal H}^{(0)}_F\right.
\simeq\bigoplus_{d\ge 0}H_F^{(d)}
\end{equation}
and
\begin{equation}
\sigma_d({\cal O}_i)=\nabla_{\delta_i}\zeta^{(-d-1)}\left|_{t'=0}
\right.,~~~i=0,\dots ,\mu-1,~d\ge 0.
\end{equation}
This is the multi-variable version of a similar formula given by Losev
\cite{lo:1}.
We can obtain the correlators with the gravitational descendants
by using the the higher residues and the above expressions.
\end{rem}
Thus we have built up the mathematical definition of the topological
Landau--Ginzburg model coupled to gravity at genus $0$:
\begin{thm}
Let $F(z,t)$ be an LG system over a frame $(Z,X,S,T)$.
Then a primitive form $\zeta^{(0)}$ uniquely defines the genus $0$ free energy
$\Phi^{grav}_F$ such that $\Phi^{grav}_F|_{t^\alpha_d=0,d>0}=\Phi_F$
up to quadratic terms.
\end{thm}
\begin{rem}
It is well known that (topological recursion relation at genus $0$
\cite{wi:1}\cite{wi:2})
\begin{equation}
\frac{\partial}{\partial t^i}\frac{\partial}{\partial t^j}\left(
\frac{\partial}{\partial t^l_{d}}\Phi^{grav}_F\right)
=\sum_{k=0}^{\mu-1} C_{ij}^k\frac{\partial}{\partial t^k}
\frac{\partial}{\partial t^0}\left(
\frac{\partial}{\partial t^l_{d}}\Phi^{grav}_F\right).
\end{equation}
If we restrict the above equation to the small phase space $H^{(0)}_F$
($t^i_d=0$,$d>0$), it becomes the ordinary Gauss--Manin differential equations
for period integrals expressed in a flat coordinate system \cite{sa:1}.
\end{rem}
\section{Primitive Forms and Mirror Symmetry}
\subsection{Period of Primitive Forms}
~~~
In this subsection, we will review the period mapping of primitive forms
\cite{sa:1}.
\begin{defn}
For any $\kappa\in\C$, we define the $D_S$-module
\begin{equation}
{\cal M}^{(\kappa)}:=D_S\left/I_S\right. ,
\end{equation}
where
\begin{equation}
I_S:=\sum_{\delta,~\delta'\in{\cal G}}D_SP(\delta,\delta')
+\sum_{\delta\in{\cal G}}D_SQ_\kappa(\delta),
\end{equation}
\begin{equation}
P(\delta,\delta'):=\delta\delta'-\delta\circ\delta'\delta_0-\nabla_\delta
\delta',~~~\delta,\delta'\in{\cal G},
\end{equation}
\begin{equation}
Q_\kappa(\delta):=w(\delta)\delta_0-\left(N-\kappa-1\right)\delta,~~~
\delta\in{\cal G}.
\end{equation}
\end{defn}
\begin{lem}
The module ${\cal M}^{(\kappa)}$ for $\kappa\in\C$ is
a simple holonomic system.
\end{lem}
For $[P]\in{\cal M}^{(\kappa+1)}$, we can define a left $D_S$-homomorphism
\begin{equation}
{\cal M}^{(\kappa+1)}\to{\cal M}^{(\kappa)},~~~[P]\mapsto[P\delta_0],
\end{equation}
since we have
\begin{equation}
P(\delta,\delta')\delta_0=\delta_0P(\delta,\delta'),
\end{equation}
and
\begin{equation}
Q_{\kappa+1}(\delta)\delta_0=\delta_0Q_\kappa(\delta).
\end{equation}
Let $l$ be the corank of the $N-\kappa-1$ on the $\mu$-dimensional horizontal
vector space of $\nabla$ on ${\cal G}$ and let $\eta_1,\dots ,\eta_l$ and
$\xi_1,\dots ,\xi_l$ be the basis of $\ker(N-\kappa-1)$ and
${\rm coker}(N-\kappa-1)$, respectively.
Then we have an $D_S$-exact sequence
\begin{equation}
0\to{\cal O}_S[E_{\kappa+1}]\oplus\bigoplus_{i=1}^l{\cal O}_S[w(\eta_i)]\to
{\cal M}^{(\kappa+1)}\to{\cal M}^{(\kappa)}\to {\cal O}_S\oplus
\bigoplus_{i=1}^l{\cal O}_S[\xi_i]\to 0,
\end{equation}
for $r\neq\kappa+1$, and
\begin{equation}
0\to\bigoplus_{i=1}^l{\cal O}_S[w(\eta_i)]\to
{\cal M}^{(\kappa+1)}\to{\cal M}^{(\kappa)}\to {\cal O}_S\oplus
\bigoplus_{i=2}^l{\cal O}_S[\xi_i]\to 0,
\end{equation}
for $r=\kappa+1$ and $\xi_1=\delta_0$, where $E_{\kappa+1}=E-(r-\kappa-1)$.
Let us define
\begin{equation}
Sol({\cal M}^{(\kappa)}):=Hom_{D_S}({\cal M}^{(\kappa)},
{\cal O}_S),~~~\kappa\in\C.
\end{equation}
\begin{rem}
Since ${\cal M}^{(\kappa)}$ has the generator $1$, a solution
$u\in Sol({\cal M}^{(\kappa)})$ is identified with $u(1)\in{\cal O}_S$.
\end{rem}
Then we have the following lemma.
\begin{lem}
The kernel of the morphism,
\begin{equation}
Sol({\cal M}^{(\kappa+1)})\stackrel{\delta_0}{\to}Sol({\cal M}^{(\kappa)}),
\end{equation}
is spanned by canonical coordinate of exponent $\kappa+1$ $($or degree
$r-\kappa$$)$ and by a constant function $1_S$.
\end{lem}
Since $Sol({\cal M}^{(\kappa)})$ always contains the constant function,
we have
\begin{equation}
dSol({\cal M}^{(\kappa)})\simeq Sol({\cal M}^{(\kappa)})\left/\C\right. ,
\end{equation}
where $dSol({\cal M}^{(\kappa)})$ is the image of $Sol({\cal M}^{(\kappa)})$
in $\Omega_S^1$ under $d$.
Let $\gamma^0(s), \dots, \gamma^{\mu-1}(s)\in H_n(X_s,\Z)$ be a $\Z$-basis of
a horizontal family of homology defined on a simply connected domain of a
covering space of $S-{\cal D}$.
We see that $\int_{\gamma(s)}\zeta^{(k-1)}\in Sol({\cal M}^{(k)})$
for $k\in\Z$.
In particular, for $k\le 0$ the constant function $1_S$ and the periods
$\int_{\gamma^i(s)}\zeta^{(k-1)}$, $i=0,\dots ,\mu-1$ form a $\C$-basis of
$Sol({\cal M}^{(k)})$.
So there exist natural isomorphisms of local systems,
\begin{align}
\dots\stackrel{\delta_0}{\simeq}dSol({\cal M}^{(-2)})|_{S-{\cal D}}
& \stackrel{\delta_0}{\simeq}dSol({\cal M}^{(-1)})|_{S-{\cal D}}\\
& \stackrel{\delta_0}{\simeq}dSol({\cal M}^{(0)})|_{S-{\cal D}}
\simeq\bigcup_{s\in S-{\cal D}}H_n(X_s,\C),
\end{align}
\begin{prop}
There exists a non-degenerate $\C$-bilinear form,
\begin{equation}
I_\kappa:dSol({\cal M}^{(\kappa)})|_{S-{\cal D}}\times
dSol({\cal M}^{(n-\kappa)})|_{S-{\cal D}}\to \C|_{S-{\cal D}},
\end{equation}
defined by
\begin{equation}
I_\kappa:=\sum_{i=1}^\mu\delta_i\otimes w(\delta^{i*})\in{\cal M}^{(\kappa)}
\otimes{\cal M}^{(n-\kappa)},~~~\kappa\in\C
\end{equation}
where $\delta_1,\dots ,\delta_\mu$ and $\delta^{1*},\dots ,\delta^{\mu *}$ are
${\cal O}_T$-basis of ${\cal G}$ and the dual basis with respect to $J$,
respectively.
\end{prop}
From this proposition, we have natural isomorphisms of local systems,
\begin{align}
\bigcup_{s\in S-{\cal D}}H^n(X_s,\C)\simeq dSol({\cal M}^{(n)})|_{S-{\cal D}}
& \stackrel{\delta_0}{\simeq}dSol({\cal M}^{(n+1)})|_{S-{\cal D}}\\
& \stackrel{\delta_0}{\simeq}dSol({\cal M}^{(n+2)})|_{S-{\cal D}}
\stackrel{\delta_0}{\simeq}\dots .
\end{align}
\begin{thm}
Let $n$ be even. For $\gamma,~\gamma'\in H_n(X_s,\Z)$, the intersection
number of these cycles is given by,
\begin{equation}
\left<\gamma,\gamma'\right>:=(2\pi\sqrt{-1})^{-n}I_0\left(d\int_\gamma
\zeta^{(-1)},d\int_{\gamma'}\zeta^{(n-1)}\right).
\end{equation}
\end{thm}
\begin{cor}
Let $n$ be even. Then
\begin{equation}\label{eq:inter}
(2\pi\sqrt{-1})^{-n}\delta_0^n:dSol({\cal M}^{(0)})|_{S-{\cal D}}
\to dSol({\cal M}^{(n)})|_{S-{\cal D}}
\end{equation}
coincide with the linear mapping defined by the intersection form.
\end{cor}
\subsection{Landau--Ginzburg Orbifolds and Calabi--Yau Varieties}
~~~
Let $f$ be a quasi--homogeneous polynomial of $z_i$ with an isolated
singularity at base point $0$ such that
\begin{equation}
f(\lambda^{w_0}z_0,\dots ,\lambda^{w_n}z_n)=\lambda f(z_0,\dots ,z_n).
\end{equation}
Let us consider the LG system $F$ associated to $f$.
Then we define the affine variety
\begin{equation}
X_1:=\{z\in\C^{n+1}|f(z_0,\dots ,z_n)+1=0\},
\end{equation}
which is the Milnor fiber of $\varphi$ at $(1,0,\dots ,0)\in S$.
We can compactify $X_1$ to $\overline{X_1}$ so that the boundary
$V:=\overline{X_1}\backslash X_1$ is given by
\begin{equation}
V=\{[z]\in \P(w_0,\dots ,w_n)|f([z])=0\}.
\end{equation}
We have two isomorphisms of $\C$-vector spaces
\begin{equation}
H^{n}(X_1,\C)\simeq \C[z_0,\dots,z_n]\left/(\frac{\partial f}{\partial z_0},
\dots ,\frac{\partial f}{\partial z_n})\right.
\end{equation}
and
\begin{equation}\label{eq:iso}
PH^{n-1}(V,\C)\simeq \C[z_0,\dots,z_n]\left/(\frac{\partial f}{\partial z_0},
\dots ,\frac{\partial f}{\partial z_n})\left|_{int}\right.\right. ,
\end{equation}
where we denote by $PH^{n-1}(V,\C)$ the primitive part of the cohomology and
$|_{int}$ means the restriction to the integral exponents sector (exponents of
the monomial $z_0^{a_0}\dots z_n^{a_n}$ is given by $\sum_{i=0}^n(a_i+1)w_i$).

\begin{rem}
Note that the integral exponents sector corresponds to the
{\it untwisted sector} of the LG orbifold theory $f//G$,
where $G$ is the discrete group acting on $z_i$ as
$z_i\to\exp(2\pi\sqrt{-1}w_i)z_i$ \cite{iv:1}\cite{ky:1}.
\end{rem}
\begin{lem}
Let $F$ be an LG system over a frame $(Z,X,S,T)$ and let us
assume that a primitive form $\zeta^{(0)}$ is given. Then
$(S',\eta_F|_{S'},\Phi_F|_{S'})$ is a Frobenius manifold, where
\begin{equation}
S':=\{t\in S|t^i=0~{\rm for}~N(\frac{\partial}{\partial t^i})\not\in\Z\}.
\end{equation}
\end{lem}
Let us assume that $r=\sum_{i=0}^nw_i=1$ and the hypersurface $V$ in weighted
projective space $\P(w_0,\dots ,w_n)$ is a smooth Calabi--Yau manifold.
Then the following conjecture is quite natural from the viewpoint of the
Landau--Ginzburg orbifold theory and the isomorphism \eqref{eq:iso}.
\begin{conj}
Let $F$ be an LG system over a frame $(Z,X,S,T)$ and let us
assume that a primitive form $\zeta^{(0)}$ is given.
Then $(S',\eta_F|_{S'},\Phi_F|_{S'})$ should be a Frobenius submanifold of the
one constructed by Barannikov--Kontsevich for a Calabi--Yau manifold $V$.
\end{conj}

Then we can formulate the mirror conjecture as follows.
\begin{conj}
There should exist a Calabi--Yau manifold $V^*$ such that $S'$ is
a Frobenius submanifold of the one constructed by Gromov--Witten theory.
\end{conj}
Indeed, if $\dim V=1$, these conjectures are true since the potential is
trivial.
\begin{rem}
Since the kernel of the map $(2\pi\sqrt{-1})^{-n}\delta_0^n$ \eqref{eq:inter}
are spanned by the canonical coordinates with integral exponents, there exist
the horizontal family of homology $\gamma^0(t),\gamma^i(t)\in H_n(X_t,\Z)$
such that
\begin{equation}
1=\int_{\gamma^0(t)}\zeta^{(0)},
\end{equation}
and
\begin{equation}
t^i=\int_{\gamma^i(t)}\zeta^{(0)},~~~{\rm for}~deg(t^i)=0.
\end{equation}
Recall that, when we calculate the `special coordinates' in the context of
mirror symmetry of Calabi--Yau manifolds, we first take the holomorphic form
$\Omega$ and calculate the period integrals by Picard--Fuchs or GKZ equations.
We choose the special solution $\int_{\Gamma^0}\Omega$ which is
single-valued near the {\it large radius limit} or {\it maximally unipotent
monodromy point} and the solutions $\int_{\Gamma^i}\Omega$ which have
simple logarithmic singularity at the maximally unipotent monodromy point.
Then we define the normalized holomorphic form by
\begin{equation}
\frac{\Omega}{\int_{\Gamma^0}\Omega},
\end{equation}
and define the `special coordinates' as
\begin{equation}
\tau^i=\frac{\int_{\Gamma^i}\Omega}{\int_{\Gamma^0}\Omega}.
\end{equation}
In \cite{bcov:1}, it is proved that the canonical coordinates at the
maximal unipotent monodromy point coincides with the `special coordinates'.
So we can conclude that to choose the holomorphic form
$\frac{\Omega}{\int_{\Gamma^0}\Omega}$ is equivalent to choose the primitive
form $\zeta^{(0)}$ at this point.
\end{rem}
\subsection{Mirror Symmetry for $\C\P^1$}
~~~
In this subsection, we construct the mirror partner of $\C\P^1$ in terms of the
primitive form.
In \cite{ba:1}, Batyrev showed that the quantum cohomology ring of $\C\P^1$ is
given by the Jacobian ring
\begin{equation}
QH^*_q(\C\P^1,\C)\simeq\C[z,z^{-1}]\left/(z\frac{\partial f}{\partial z})
\right. ,
\end{equation}
where
\begin{equation}
f:=z+qz^{-1},~~~0<|q|<1.
\end{equation}
This fact implies that it may be possible to construct a primitive form
associated to the holomorphic function
\begin{equation}
f:\C^*\to \C.
\end{equation}
Let $(z,t^1)$, $(t^0,t^1)$ and $(t^1)$ be the coordinate of
$X':=\C^*\times\C$, $S':=\C^2$ and $T':=\C$, respectively.
Then we define the function on $Z':=X'\times_{T'}S'$
\begin{equation}
F(z,t):=t^0+z+q\exp(t^1)z^{-1},
\end{equation}
where $q$ is a parameter of a punctured unit disk $0<|q|<1$.
Then we can define a frame
\begin{equation}
\begin{CD}
Z' @>\hat{\pi'}>> X' \\
@V p' VV @VV q' V \\
S' @>\pi'>> T' ,
\end{CD}
\end{equation}
the map $\varphi':X'\to S'$, the critical set $C\subset X'$ and the
discriminant $D=\varphi'(C)$.
In the previous sections, we have assumed that $f$ is a quasi--homogeneous
polynomial of $z_i$ with an isolated singularity at base point $0\in\C^{n+1}$.
In such case, the maps $\pi$ and $q$ are Stein, which are trivial fibrations
with contractible fibers and the map $\varphi$ is Stein.
This condition is sufficient for the existence of higher residue pairings
$K^{(k)}$.
The map $q'$ and $\varphi'$ do not satisfy the above conditions, but it is
possible to replace a frame $(Z',X',S',T')$ by $(Z,X,S,T)$ such that
this new frame satisfies the conditions.
Indeed, we can take a new frame $(Z,X,S,T)$ as follows.
\begin{equation}
T:=\{(t^1)\in T':||(t^1)-(0)||<\epsilon\},
\end{equation}
\begin{equation}
\begin{split}
S:= & \{(t^0,t^1)\in S':||(t^0,t^1)-(2q^{\frac{1}{2}},0)||<\delta_1\}\\
& \cup\{(t^0,t^1)\in S':||(t^0,t^1)-(-2q^{\frac{1}{2}},0)||<\delta_2\}
\cap{\pi'}^{-1}(T),
\end{split}
\end{equation}
\begin{equation}
\begin{split}
X:= & \{(z,t^1)\in X':||(z,t^1)-(q^{\frac{1}{2}},0)||<r_1\}\\
& \cup\{(z,t^1)\in X':||(z,t^1)-(q^{\frac{1}{2}},0)||<r_2\}
\cap{\varphi'}^{-1}(S),
\end{split}
\end{equation}
\begin{equation}
Z:=X\times_T S,
\end{equation}
where $\epsilon$, $\delta_i$ and $r_i$, $(i=1,2)$ are small real number
with $1>>r_i>>\delta_i>>\epsilon>0$ (i=1,2) and $||\cdot ||$ is some Euclidean
distances.
The maps are given by restriction .
We see easily that the correspondence
\begin{equation}
T(S)_0\simeq \C\frac{\partial}{\partial t^0}\oplus
\C\frac{\partial}{\partial t^1}\to
\C[z,z^{-1}]\left/(z\frac{\partial f}{\partial z})\right.,~~~
\delta\mapsto\hat{\delta}F|_{p^{-1}(0)},
\end{equation}
is bijection, so we may call $F$ an LG system.
\begin{thm}
There exist higher residue pairing $K^{(k)}$ such that
\begin{equation}
\begin{aligned}
K^{(0)}([\phi_1\frac{dz}{z}],[\phi_2\frac{dz}{z}]) & = \sum_{z\frac{\partial F}
{\partial z}=0}\frac{\phi_1\phi_2}
{z\frac{\partial}{\partial z}(z\frac{\partial F}{\partial z})}\\
&=\frac{1}{2\pi\sqrt{-1}}\int_\gamma\frac{\phi_1\phi_2}
{z\frac{\partial F}{\partial z}}\frac{dz}{z},
\end{aligned}
\end{equation}
for $[\phi_1dz/z]$, $[\phi_2dz/z]\in\pi_*{\cal H}^{(0)}_F$ and $\gamma$ is a
homology cycle of $X\backslash\{ z\frac{\partial F}{\partial z}=0\}$, defined
by a general fiber of the mapping
\begin{equation}
X\to \R^{+}\times T,~~~(z,t^1)\mapsto (|z\frac{\partial F}{\partial z}|,t^1).
\end{equation}
\end{thm}
\begin{thm}
The primitive form for $F(z,t)=t^0+z+q\exp(t^1)z^{-1}$ is
\begin{equation}
\zeta^{(0)}=\frac{dz}{z}.
\end{equation}
and the minimal exponent is $r=0$, the Euler vector field is
\begin{equation}
E=t^0\frac{\partial}{\partial t^0}+2\frac{\partial}{\partial t^1},
\end{equation}
and the canonical coordinates are given by $(t^0,t^1)$.
\end{thm}
\begin{pf}
We will show that $[\frac{dz}{z}]$ satisfy the five conditions of Definition 1.4.
\begin{enumerate}
\item It is obvious.
\item Since $K^{(1)}$ is skew symmetric, we only have to show that
\begin{equation}
K^{(1)}(\nabla_{\delta_0}\zeta^{(-1)},\nabla_{\delta_1}\zeta^{(-1)})=0.
\end{equation}
This follows from
\begin{equation}
\begin{aligned}
K^{(1)}(\nabla_{\delta_0}\zeta^{(-1)},\nabla_{\delta_1}\zeta^{(-1)}) &=
\frac{1}{2\pi\sqrt{-1}}\cdot(-\frac{1}{2})\int_\gamma\frac{q\exp(t^1)z^{-1}}
{(z\frac{\partial F}{\partial z})^2}\frac{dz}{z}\\
& = 0.
\end{aligned}
\end{equation}
\item The Euler vector field $E$ is given by
\begin{equation}
\begin{aligned}
E& =t^0\frac{\partial}{\partial t^0}-t^0\circ\frac{\partial}{\partial t^0}\\
& =t^0\frac{\partial}{\partial t^0}+2\frac{\partial}{\partial t^1}.
\end{aligned}
\end{equation}
So we have
\begin{equation}
\begin{aligned}
\nabla_E[\frac{dz}{z}]& =t^0\nabla_{\delta_0}[\frac{dz}{z}]+2\nabla_{\delta_1}
[\frac{dz}{z}]\\
&=-\frac{(z+q\exp(t^1)z^{-1})^2}{(z-q\exp(t^1)z^{-1})^2}\frac{dz}{z}
 +2\frac{2q\exp(t^1)}{(z-q\exp(t^1)z^{-1})^2}\frac{dz}{z}\\
&=-\frac{dz}{z}.
\end{aligned}
\end{equation}
This shows that $r=0$.
\item Since
\begin{equation}
\frac{\partial }{\partial t^1}\circ\frac{\partial}{\partial t^1}=
q\exp(t^1)\frac{\partial}{\partial t^0 }
\end{equation}
and
\begin{equation}
\begin{aligned}
z\frac{\partial}{\partial z}\left(\frac{q\exp(t^1)z^{-1}\cdot q\exp(t^1)z^{-1}
-q\exp(t^1)}{z\frac{\partial F}{\partial z}}\right)& = q\exp(t^1)z^{-1}\\
& = \frac{\partial }{\partial t^1}(\frac{\partial F}{\partial t^1}),
\end{aligned}
\end{equation}
we have
\begin{equation}
\begin{aligned}
\nabla_{\delta_1}&\nabla_{\delta_1}\nabla_{\delta_0}^{-1}[\frac{dz}{z}] -
\nabla_{\delta_1\circ\delta_1}[\frac{dz}{z}]\\
&= \nabla_{\delta_0}\left(
q\exp(t^1)z^{-1}\cdot q\exp(t^1)z^{-1}-q\exp(t^1)\frac{dz}{z}\right) +
\frac{\partial }{\partial t^1}(\frac{\partial F}{\partial t^1})\frac{dz}{z}\\
&= -z\frac{\partial}{\partial z}\left(\frac{q\exp(t^1)z^{-1}\cdot q\exp(t^1)
z^{-1}-q\exp(t^1)}{z\frac{\partial F}{\partial z}}\right)\frac{dz}{z}+
\frac{\partial }{\partial t^1}(\frac{\partial F}{\partial t^1})\frac{dz}{z}\\
&= 0.
\end{aligned}
\end{equation}
It is clear that
\begin{equation}
\nabla_{\delta_0}\nabla_{\delta_1}\nabla_{\delta_0}^{-1}[\frac{dz}{z}] -
\nabla_{\delta_0\circ\delta_1}[\frac{dz}{z}]=0
\end{equation}
and
\begin{equation}
\nabla_{\delta_0}\nabla_{\delta_0}\nabla_{\delta_0}^{-1}[\frac{dz}{z}] -
\nabla_{\delta_0\circ\delta_0}[\frac{dz}{z}]=0.
\end{equation}
\item Since we have proved that $K^{(1)}(\nabla_\delta\zeta^{(-1)},
\nabla_{\delta'}\zeta^{(-1)})=0$, the condition becomes the following form;
\begin{equation}
t^0\nabla_\delta[\frac{dz}{z}]=\nabla_{t^0\circ\delta}[\frac{dz}{z}]+
\nabla_{(N-1)\delta}\nabla_{\delta_0}^{-1}[\frac{dz}{z}],~~~\delta\in {\cal G}.
\end{equation}
We have
\begin{equation}
t^0\circ\frac{\partial}{\partial t^1}=
-2q\exp(t^1)\frac{\partial}{\partial t^0}.
\end{equation}
Since
\begin{equation}
t^0\nabla_{\delta_1}[\frac{dz}{z}]
=-\frac{2q\exp(t^1)(z+q\exp(t^1)z^{-1})^2}{(z-q\exp(t^1)z^{-1})^2}
\frac{dz}{z}
\end{equation}
and
\begin{equation}
2q\exp(t^1)\nabla_{\delta_0}[\frac{dz}{z}]=
+\frac{2q\exp(t^1)(z+q\exp(t^1)z^{-1})^2}{(z-q\exp(t^1)z^{-1})^2}\frac{dz}{z},
\end{equation}
we have
\begin{equation}
t^0\nabla_{\delta_1}[\frac{dz}{z}]=\nabla_{t^0\circ\delta_1}[\frac{dz}{z}].
\end{equation}
\end{enumerate}
Now, it is clear that the $(t^0,t^1)$ is a canonical coordinate system.
\qed
\end{pf}
We can easily show that the potential $\Phi_F$ is given by
\begin{equation}
\Phi_F=(t^0)^2t^1+q\exp(t^1)
\end{equation}
up to quadratic terms.
So we have constructed the mirror partner of $\C\P^1$.
\begin{thm}
Let $f(z)=z+qz^{-1}$ be a holomorphic function on $\C^*$.
Then the pair $(\C^*,f)$ is a mirror manifold of $\C\P^1$.
In other words, we have the equation
\begin{equation}
\Phi_{\C\P^1}^{st}(t)=\Phi_F^{grav}(\tilde{t}),
\end{equation}
where $\Phi_{\C\P^1}^{st}$ is the $($large phase space$)$ stable genus $0$
generating function for Gromov--Witten Invariants
\begin{equation}
\Phi_{\C\P^1}^{st}(t):=\sum_{n\ge 3,(a_i,d_i)}\frac{1}{n!}t^{a_1}_{d_1}\dots
t^{a_n}_{d_n}\sum_{\beta\in \Z_{\ge 0}}q^\beta
\left<\sigma_{d_1}({\cal O}_{a_1})\dots\sigma_{d_n}({\cal O}_{a_n})
\right>_{0,\beta[\C\P^1]}
\end{equation}
and
\begin{equation}
\tilde{t}^i_d:=t^i_d+\sum_{(j,e)}t^j_e\sum_{\beta\in\Z_{\ge 0}}q^\beta
\left<\sigma_{e-d-1}({\cal O}_j){\cal O}^i\right>_{0,\beta[\C\P^1]}.
\end{equation}
\end{thm}
The quantum cohomology ring of a toric Fano manifold $\P_\Sigma$ associated to
the $n$-dimensional fan $\Sigma$ is given by the Jacobian ring \cite{ba:1}
\begin{equation}
QH^*_\varphi(\P_\Sigma,\C)\simeq\C[z_1^{\pm},\dots ,z_n^{\pm}]\left/
(z_1\frac{\partial f}{\partial z_1},\dots ,
z_n\frac{\partial f}{\partial z_n})\right. ,
\end{equation}
of the Laurent polynomial
\begin{equation}
f_\varphi(z)=\sum_{i=1}^d\exp(\varphi(v_i))^{-1}z^{v_i},
\end{equation}
where ${v_1,\dots ,v_d}$ is the set of all generators of $1$-dimensional cones
in $\Sigma$ and $\varphi$ is an element of the complexified K\'{a}hler cone
of $\P_\Sigma$.
It may be possible to construct a primitive form in these cases, but we are
confronted by two difficulties which do not exist in $\C\P^1$ case.
The first is to find a frame $(Z,X,S,T)$ and the second is to find a (good) LG
system $F(z,t)$.
However, from the results \cite{ehx:2}, we can imagine that
\begin{conj}
There should exist a primitive form for $f$ which is given by
$\zeta^{(0)}=[\frac{dz_1}{z_1}\wedge\dots\wedge\frac{dz_n}{z_n}]$.
Then the pair $((\C^*)^n,f)$ is a mirror manifold of a toric Fano manifold
$\P_\Sigma$.
\end{conj}
It will be important to use the notion of primitive forms if one consider the
more general mirror problems.
For example, ADE-type theory does not correspond to the usual manifolds,
but there may exist a `mirror symmetry' \cite{ky:1}\cite{ta:1} and a structure
of quantum cohomology such as Virasoro conditions \cite{ehx:1}.
%

%
\end{document}